\documentclass{article} 

\usepackage{amsmath}
\usepackage{amsthm}
\usepackage{graphicx}
\usepackage{float}
\newtheorem{prop}{Proposition}

\newtheorem{cor}{Corollary}

\usepackage{xcolor}

\usepackage[justification=centering]{caption}
\usepackage{subcaption}

\begin{document}

\author{Elizabeth J. Schaefer \thanks{Yale University elizabeth.schaefer@yale.edu} \\ Andrew J. Schaefer \thanks{Rice University andrew.schaefer@rice.edu}}
%
\title{Maximizing the  Score in ``Ticket to Ride''}
%
%
%

%
%

\markboth{IEEE TRANSACTIONS ON GAMES}%
{Shell \MakeLowercase{\textit{et al.}}: Bare Demo of IEEEtran.cls for IEEE Journals}
%



\maketitle

\begin{abstract}
We give two graph-theoretic models and a mixed-integer program to calculate the maximum achievable score in the popular board game ``Ticket to Ride.'' In Ticket to Ride, players compete to claim railway routes on a map, with points awarded based on the length of each route and the successful completion of destination tickets connecting specific city pairs. Each player has 45 train cars available, and each route can be chosen by only one player. Using the mixed-integer programming model, we examine the optimal solution with the 45 allocatable train cars, leading to an optimal score of 285 points. We also calculate the optimal solutions for up to 50 train cars. We determine the most frequently chosen tickets and routes over these 50 instances, giving insight into how optimization might be used to balance games. In particular, we identify several instances in which the point values can be adjusted to better balance the game.
\end{abstract}

\vspace{-0.6em}

%

\section{Introduction}
%
%
%
%
``Ticket to Ride,'' by Days of Wonder, is a popular board game set on numerous countries and continents. We focus on the original version, which is set on a North American railroad network. The game won the 2004 Meeples’ Choice Award, the 2005 International Gamers Award for General Strategy - Multiplayer, the 2004 Spiel des Jahres, the Origins Award for Best Board Game of 2004, and the 2005 Diana Jones award\footnote{https://www.originsawards.net/origins-award-winners}. Over eight million copies of Ticket to Ride and its variants have been sold\footnote{http://www.internationalgamersawards.net/}.  Players compete to establish train routes between various destinations, and points are earned based on a combination of the length of the claimed routes, the longest continuous railway, and whether the player completed a connection between various city tickets. 

The objective of the game is to score the highest number of points. Each turn, a player can do one of three things: 1) draw one or two random train car cards, which are the ``currency'' to claim a route;  2) claim an available route by discarding an appropriate set of train cards and placing the requisite number of train cars on it, thus earning points that correspond to the length of the route and helps complete destination tickets; or 
3) draw three additional tickets, keeping at least one of them. Destination tickets are a \textcolor{black}{unique} pair of cities, and if the player \textcolor{black}{successfully} completes a passage between these cities, \textcolor{black}{no matter the length}, she receives the associated point value. However, if the player fails to complete a route, its point value is deducted from her score. The game ends when one player has two or fewer train cars left, and all players get one final turn, and the player with the highest score wins. In the North American version, the player with the longest continuous railway is awarded an extra 10 points\footnote{The full rules are available at https://ncdn0.daysofwonder.com/tickettoride/ en/img/7201-T2R-Rules-EN-2019.pdf}. In this paper we focus only on the points earned by claiming routes and satisfying tickets.  

This paper considers the question of what is the maximum possible score in Ticket to Ride. As such, we ignore competition,  as players can only reduce each other's scores by competing for routes and tickets. Most routes and all tickets may be claimed by only one player. 

Others have modeled Ticket to Ride. Witter and Lyford \cite{TealWitter20} proposed a method for assigning scores to tickets using probability and graph theory. de Mesentier Silva et al. \cite{Mesentier17} considered four agents with different playing strategies, and compared their outcomes. They then used an evolutionary algorithm to generate boards and tickets \cite{Mesentier18}. Yang et al. apply reinforcement learning to a scaled-down version of Ticket to Ride, and show that RL-trained agents can beat heuristic strategies \cite{Yang}.


\section{Graph Representation and Optimization Model}   

\subsection {Undirected Graph Representation}
We define the Ticket to Ride problem over a set $V$ of vertices and a set $E$ of edges that define an undirected \textcolor{black}{connected} graph $G = (V,E)$. We assume that $G$ contains no self-loops. Each edge $e \in E$ has an associated length $l_e > 0$ and point value $\color{black}{f_e} > 0$. In the board game Ticket to Ride, there is a single function mapping the edge lengths $l$ to the point values $c$, \textcolor{black}{(e.g., every edge of length two is worth two points),} but our model does not make this assumption. For convenience we assume a linear ordering on the vertices and edges, and for each edge $e=\color{black}{\{i,j\}}$,  $i < j$. There is a set $K$ of ``tickets,'' where each ticket $k \in K$ is defined by a pair of nodes $s_k \in V$ and $t_k \in V$, where $s_k < t_k$, and an associated point value $d_k > 0$. We assume that for each ticket $k$, $s_k$ and $t_k$ are in the same connected component of $G$; otherwise the ticket may be eliminated as it is impossible to satisfy.  \textcolor{black} {We assume that no ticket coincides with a single edge, which is the case in the board game Ticket to Ride. We also assume that all routes that could be claimed by two different players are simplified to a single edge.}
 
Consider a subset of edges $E' \subseteq E$, where the total length of all edges in $E'$ does not exceed the budget $\beta > 0$, that is, $\sum_{e \in E'} l_e \ \leq \ \beta$. \textcolor{black}{In the board game Ticket to Ride this captures the 45 train cars allocated to each player.} Define the (possibly disconnected) subgraph edge-induced by $E'$ as $G' = (V',E')$, where $V' \subseteq V$ consists of those vertices adjacent to at least one edge in $E'$. This in turn defines a subset of ``completed tickets'' $\begin{color}{black}K(E') \end{color} \subseteq K$ defined as 
\[
\begin{color}{black}K(E') \end{color}  \ = \ \{k \in K \ | \ s_k \ \mathrm{and} \ t_k \mathrm{ \ are \ connected \ in \ } G' \}. 
\]
Given a subset of edges $E'$, its ``total points'' $C(E')$ is the sum of the point values of each edge in $E'$ plus the sum of the point values over $\begin{color}{black}K(E') \end{color}$, the subset of every completed ticket in the edge-induced subgraph $G'$. That is,
\[
\begin{color}{black}F(E')\end{color} \ = \ \sum_{e \in E'} \begin{color}{black}f_e \end{color} \ + \ \sum_{k \in \begin{color}{black}K(E') \end{color}} d_k.
\]
 
The Ticket To Ride problem is to choose a subset of edges  whose total length does not exceed $\beta$ with the highest total points. That is, to solve 
\[
\max_{E' \subseteq E} \left \{ \begin{color}{black}F(E') \end{color} \ \left | \ \sum_{e \in E'} l_e \ \leq \ \beta \right \} \right.. 
\]

\begin{figure}[h!]
    \centering
    \begin{subfigure}[t]{.22\textwidth}
        \centering
        \includegraphics[scale=.12]{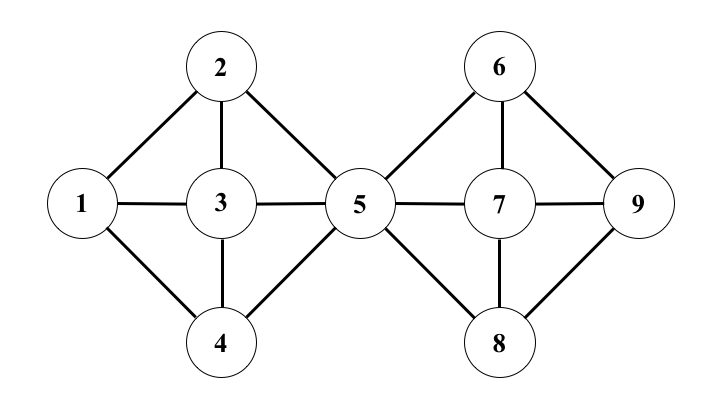}
        \caption{An example of the undirected graph formulation of the Ticket to Ride problem.}
        \label{UndirectedGraphExample}
    \end{subfigure}%
    \hspace{5pt}
    \begin{subfigure}[t]{.22\textwidth}
        \centering
        \includegraphics[scale=.12]{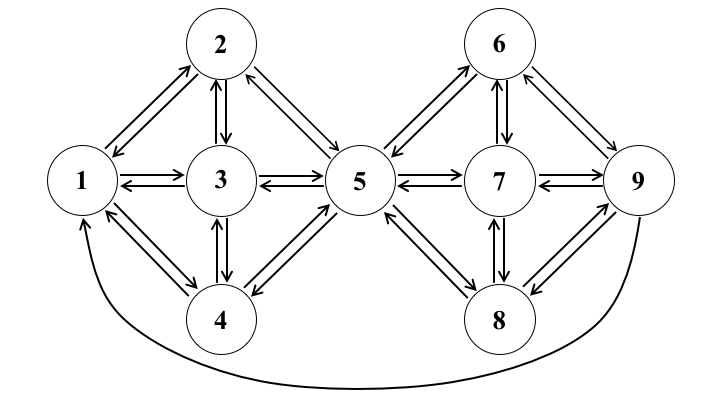}
        \caption{The directed graph formulation. There is a single ticket defined on cities 1 \& 9.}
        \label{DirectedGraphExample}
    \end{subfigure}
    \captionsetup{justification=raggedright,singlelinecheck=off}
    
    \caption{Graph formulations of the Ticket to Ride problem. Each edge in Figure \ref{UndirectedGraphExample} is replaced by a forward and backward arc, and there is a ``dummy'' arc from 9 to 1 corresponding to the ticket defined on cities 1 to 9. Figure \ref{DirectedGraphExample} shows the transformation from the undirected graph given in Figure \ref{UndirectedGraphExample} to a directed graph.}
    \label{GraphExample}
\end{figure}

Figure \ref{GraphExample} gives an illustrative example of the undirected \textcolor{black}{and directed} graph formulation; the only ticket is from City 1 to City 9. 

Formally, given an instance of the Ticket to Ride problem and an integer $\alpha$, the Ticket to Ride feasibility question is: 

\texttt{
TICKET TO RIDE:}
Does there exist a subset of edges $E' \subseteq E$ whose total length does not exceed $\beta$, and a subset of tickets $\begin{color}{black}K(E') \end{color} \subseteq K$ of tickets connected in $E'$,  whose total points $\color{black}{F(E')}$ is at least $\alpha$?  

\texttt{TICKET TO RIDE} is trivially NP-hard, as it generalizes $0$-$1$ \texttt{KNAPSACK} \cite{GareyJohnson}. 



\subsection{\textcolor{black}{Relationship to }\textcolor{black}{Other Problems in the Literature}}
\begin{color}{black} 
Ticket to Ride falls into the broad category of fixed-charge network design problems, and  multicommodity network design problems (the literature is vast; see \cite{crainic2020fixed} and \cite{gendron1999multicommodity}, respectively, for  surveys). In the former, an edge or arc must be constructed at a fixed cost before it may be used. In the latter, $K$ commodities are given as source-destination pairs, which must be connected in a network. The tickets in Ticket to Ride are naturally formulated as commodities. Unlike in standard multicommodity network design problems, in Ticket to Ride, not all commodities (tickets) must be connected, and connected tickets yield a ``prize'' (the bonus for connected tickets).  In this sense, within the broad classes of fixed-charge network design problems and multicommodity design problems, Ticket to Ride has similarities to prize-collecting Steiner tree problems (e.g., \cite{ljubic2006algorithmic}) and prize-collecting Steiner forest problems  
(e.g., \cite{Ahmadi25}). 
\end{color}

\subsection {Directed Graph Formulation} \label{DirectedGraphModel} We formulate the Ticket to Ride problem as  a mixed-integer program on a family of directed graphs $\{H_k\}_{k \in K}$.  Given ticket $k \in K$, we construct $H_k$, by replacing each undirected edge $e = \color{black}{\{i,j\}}$ (recall $i<j$ for each edge), with a ``forward'' directed arc $a = (i,j)$ and a ``backward'' directed arc $a = (j,i)$.   We also create a dummy backward arc $(t_k,s_k)$ connecting the two cities for ticket $k$ (recall that $s_k < t_k$). These dummy arcs are used for determining the completion of tickets. Let $\mathcal{A}^k$ be the set of arcs, which includes the one or two directed arcs for each original edge $e \in E$, and the dummy arc $(t_k,s_k)$. Define the directed graph $H_k \ = \ (V,\mathcal{A}^k)$. Given a subset of edges $E' \subseteq E$, let $H'_k$ be the directed subgraph of $H'$ induced by the edges in $E'$.


\begin{prop} \label{cycleprop}
Given a subset of edges $E' \subseteq E$, for every ticket $k \in K$, there exists a directed cycle in $H'_k$ that contains the arc $(t_k,s_k)$ if and only if nodes $s_k$ and $t_k$ are connected in the undirected subgraph $G'$.
\end{prop}

\begin{proof}
``$\Rightarrow$'':  Given a directed cycle in $H'_k$  containing  arc $(t_k,s_k)$, delete  arc $(t_k,s_k)$. Edges corresponding to  remaining arcs connect $s_k$ and $t_k$ in the undirected subgraph $G'$.

``$\Leftarrow$'': If nodes $s_k$ and $t_k$ are connected in the undirected subgraph $G'$, let $P$ be an $s_k-t_k$ path in $G'$. The directed $(s_k,t_k)$ path in $H'_k$ corresponding to $P$ plus the directed arc $(t_k,s_k)$ forms a directed cycle in $H'_k$. 
\end{proof}

\begin{cor} \label{CycleExistenceCorollary} Given a subset of edges $E' \subseteq E$, the set of completed tickets $\begin{color}{black}K(E') \end{color} \subseteq K$ can be expressed as:
\begin{gather*}
\begin{color}{black}K(E') \end{color} = \{k \in K \mid \exists \ \mathrm{a\ directed\ cycle\ in\ } H'_k \\
\mathrm{containing\ arc\ } (t_k,s_k)\}.
\end{gather*}
\end{cor}  

\section{Mixed-Integer Programming Model}
We model the directed graph formulation of the Ticket to Ride problem given in Section \ref{DirectedGraphModel} as a mixed-integer program. A  linear mixed-integer program is an optimization problem that maximizes (or minimizes) a linear function of decision variables subject to a set of linear constraints. Some of the decision variables must take on integer values, while the others are continuous. There are deep connections between graph theory and mixed-integer programming; see Wolsey \cite{Wolsey} for an introduction.


The decision variables are:
\begin{itemize}

\item $x_e \ = \ \left \{ 
\begin{tabular}{ ll } 
 1 & if edge  $e$  is chosen, \\ 
 0 & otherwise.
\end{tabular}
\right.$

\item $y_{(i,j),k} \ = $ the amount of flow  sent along arc $(i,j)$ that is used to satisfy ticket $k$.
\end{itemize}
Notice that the $x$ variables are defined with respect to the undirected graph $G$, as flow for a given edge $e$ may be used for different directions for different tickets.

Consider the following mixed-integer program:

\begin{subequations} \label{IPmodel}
 \begin{equation}  \label{IPOjective}
     \max  \sum_{e \in E} \begin{color}{black}f_e \end{color} x_e \ + \ \sum_{k \in K} d_k y_{(t_k,s_k),k} 
 \end{equation}
 \noindent subject to
 \begin{equation}  \label{IPLength}
     \sum_{e \in E} l_e x_e \ \leq \ \beta,
 \end{equation}
 \begin{equation} 
 \label{IPVarUB}
   y_{(i,j),k}  \ + \  y_{(j,i),k}  \ \leq \ x_e, \ \mathrm{for} \ e = (i,j) \in E, 
 \end{equation}
 
\begin{equation}  \label{IPFlowBalance}
     \sum_{\substack{(j,i) \in \mathcal{A}^k}} y_{(j,i),k} - \sum_{\substack{(i,j) \in \mathcal{A}^k}} y_{(i,j),k} = 0, \forall \ k \in K, \forall \ i \in V,
 \end{equation}
\begin{equation}  
\label{IPSignRestrictions}
 x \in \{0,1\}^{|E|}, \ y_{(i,j),k} \in \{0,1\}, \forall k \in K, \forall (i,j) \in \mathcal{A}^k.
 \end{equation} 
\end{subequations}

Objective (\ref{IPOjective}) maximizes the total number of points earned by claiming routes plus those earned by satisfying tickets.  Constraint (\ref{IPLength}) ensures that the sum of the lengths of the chosen routes may not exceed $\beta$.  Constraints (\ref{IPVarUB}) allow flow on an (non-dummy) arc only if the corresponding edge is chosen, and, if an edge is chosen, flow is permitted in at most one direction. Notice that dummy arcs do not appear in  (\ref{IPVarUB}). Constraints (\ref{IPFlowBalance}) ensure that for every ticket $k \in K$, the flow entering node $i$ must equal the flow leaving node $i$. Recall that the ``dummy arcs'' $(t_k,s_k)$ belong to the set $\mathcal{A}^k$. It can be shown that integrality on the $y$ variables can be relaxed.

\begin{prop}
Given an optimal solution $(x^*,y^*)$ to the mixed-integer program (\ref{IPmodel}), an optimal solution to the Ticket to Ride problem is to choose edges  $E' \ = \ \{e \in E \ | \ x^*_e \ = \ 1\}$, thereby satisfying tickets $\begin{color}{black}K(E') \end{color} \ = \ \{k \in K \ | \ y^*_{(t_k,s_k),k} \ = 1$\} and giving $c^T x^* \ + \ d^T y^* \ = \ F (E')$ points.
\end{prop}

\begin{proof}
Clearly the solution described in the proposition is feasible for the Ticket to Ride problem, as  $x^*$  satisfy the budget constraints and each ticket  $k \in \begin{color}{black}K(E') \end{color}$ must be satisfied, as $y^*_{(t_k,s_k),k} \ = 1$ implies that $s_k$ and $t_k$ are connected in $E'$. Suppose there exists a solution $\overline{E},\overline{K}$ to the Ticket to Ride problem that gives strictly more points than $F (E'))$. Then consider the solution $\overline{x}_j \ = \ 1$ for $j \in \overline{E},$ $0$ otherwise, and $\overline{y}_{(a,k)} = 1$ for every $k \in \overline{K}$ and all arcs $a$ on a cycle formed by the dummy arc $(t_k,s_k)$ and any directed path in $\overline{E}$ from $s_k$ to $t_k$ (which must exist, as ticket $k$ is satisfied). This solution $\overline{x}, \overline{y}$ is feasible for the mixed-integer program (\ref{IPmodel}) and gives a higher objective than the optimal solution, a contradiction.
\end{proof}

\begin{cor}
The optimal objective of the mixed-integer program (\ref{IPmodel}), plus the ten-point bonus for the longest railway, provides an upper bound on the score of the competitive version of the board game Ticket to Ride. 
\end{cor}

\textcolor{black}{Finally, we note that our model could be used to maximize the possible achievable score of a partially completed game, even in the presence of competition. At any given point in a partially completed game, modifying  model (\ref{IPmodel}) with the remaining number of train cars $\beta' \leq \beta$,  the set of unclaimed routes $K' \subseteq K$, and the still-available   edges $E' \subseteq E$ will yield the maximum possible remaining score in the partially completed game.}

\section{Computational Study}
The Ticket to Ride game consists of 36 cities, 78 routes, and 30 tickets.   Model (\ref{IPmodel}) has 78 binary decision variables, 4710 continuous decision variables, and 1,159 constraints. We used Gurobi 9.1, a commercial MIP solver\footnote{www.gurobi.com}.

\subsection{The Optimal Solution for 45 Train Cars}
\begin{figure}
\centering
\includegraphics[scale=.30]{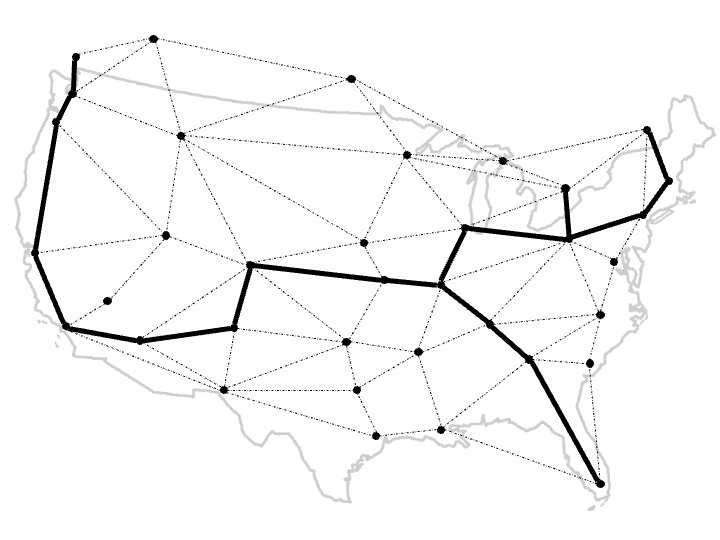}
\caption{Map of the optimal solution with 45 train cars.} 
\label{OptimalMap}
\end{figure}

The rules of Ticket to Ride give each player 45 train cars to be allocated across different routes. The optimal score is 285, \textcolor{black}{of which 223 points comes from completed tickets, and}  ignoring the ten-point bonus for longest railway (see Figure \ref{OptimalMap}). \textcolor{black}{A set of optimal routes are shown in Table \ref{OptimalRoutes} and the optimal tickets are shown in Table \ref{ReduceTicketValue}. Notice that the average route length is 2.5. These solutions reflect the solution provided by Gurobi on default settings; alternative optimal solutions were not considered.}

\begin{table}
\color{black}
 \centering 
    \begin{tabular}{|l|c|}
    \hline
    Route &  Length \\
    \hline
    Vancouver to Seattle & 1 \\
    Seattle to Portland & 1 \\
    Portland to San Francisco & 5 \\
    San Francisco to Los Angeles & 3 \\
    Los Angeles to Phoenix & 3 \\
    Phoenix to Santa Fe & 3 \\
    Santa Fe to Denver & 2 \\
    Denver to Kansas City & 4 \\
    Kansas City to St. Louis & 2 \\
    St. Louis to Nashville & 2 \\
    Nashville to Atlanta & 1 \\
    Atlanta to Miami & 5 \\
    St. Louis to Chicago & 2 \\
    Chicago to Pittsburgh & 3 \\
    Pittsburgh to Toronto & 2 \\
    Pittsburgh to New York & 2 \\
    New York to Boston & 2 \\
    Boston to Montreal & 2 \\
    \hline
 \end{tabular}
    \caption{\color{black}The \textcolor{black}{set of routes in an optimal solution} for 45 train cars.} \label{OptimalRoutes}   
\end{table}

\subsection{Sensitivity Analyses}
We first performed a sensitivity analysis on the ticket values of the North American 45-train-car version. \begin{color}{black} Table \ref{ReduceTicketValue} displays the change in the optimal objective if a chosen ticket were removed.   \end{color}


\begin{table}[]
\color{black}
    \centering 
    \begin{tabular}{|l|c|c|}
    \hline
    &  & Impact of \\
    Ticket    &  Points & Ticket Removal\\
    \hline
    Atlanta to New York & 6 & 0\\
    Atlanta to Montreal & 9  & - 9\\
    Los Angeles to Seattle & 9  & -9\\
    Chicago to Santa Fe & 9  & -9\\
    Miami to Toronto & 10  & - 1\\
    Phoenix to Portland & 11  & -5\\
    Denver to Pittsburgh & 11  & 0\\
    Boston to Miami & 12  & - 1\\
    Santa Fe to Vancouver & 13  & - 13\\
    Chicago to Los Angeles & 16 & -15\\
    Atlanta to San Francisco & 17 & -13\\
    Nashville to Portland & 17  & 0\\
    Los Angeles to Miami & 20  & -1\\
    Montreal to Vancouver & 20  & -20\\
    Los Angeles to New York & 21 & -21\\
    New York to Seattle & 22  & -1\\
    \hline
 \end{tabular}
    \caption{\color{black} The \textcolor{black}{set of} tickets  \textcolor{black}{in an optimal solution} for 45 train cars, and impact on the optimal score if each ticket were removed.} \label{ReduceTicketValue}
\end{table}


\begin{figure}
\centering

\includegraphics[scale=.25]{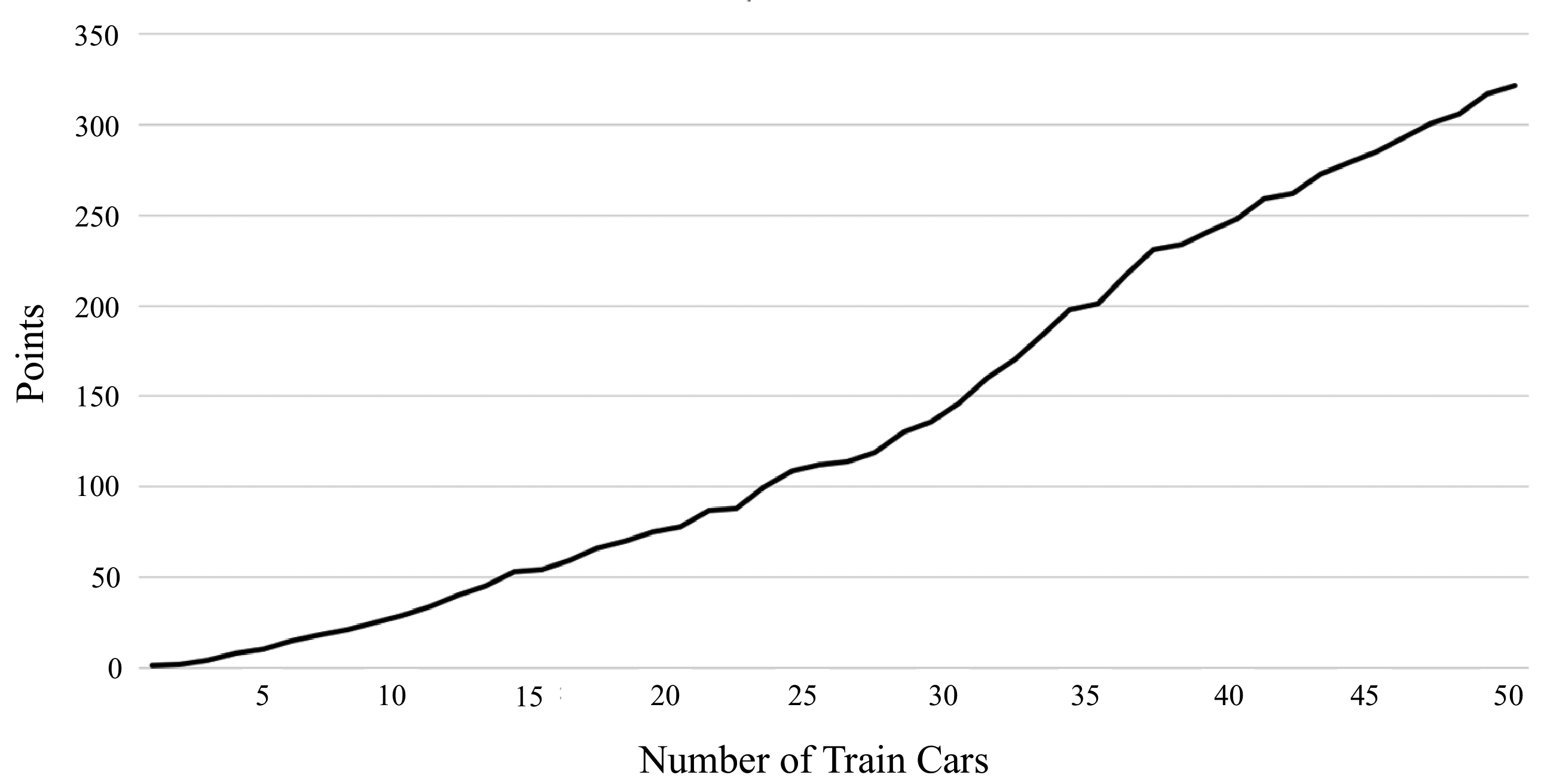}
\caption{Optimal scores as the number of train cars varies.} 
\label{SummaryFigure}
\end{figure}

We solved the model (\ref{IPmodel}) for the number of train cars $\beta$ ranging from 1 to 50; recall that the game has $\beta = 45$ (see 
Figure \ref{SummaryFigure}). We then determined which routes and tickets were chosen most frequently in \textcolor{black}{ the 50 optimal solutions (ignoring alternate optima)}. Table \ref{CardFrequency} shows the ten most frequently selected tickets. Interestingly, all ten are selected in the optimal solution with 45 train cards. While the average point value of all 30 tickets is 11.63, the average point value of the ten most frequently chosen tickets is 14.1. There are seven tickets with point values of 16 or higher; five of these seven are among the most frequently chosen, indicating that a more balanced game might result from reducing the point values of these tickets. It also provides evidence that choosing high-value tickets may be a good strategy, although this insight may not extend to a competitive environment.

\begin{figure}
\centering
\includegraphics[scale=.30]{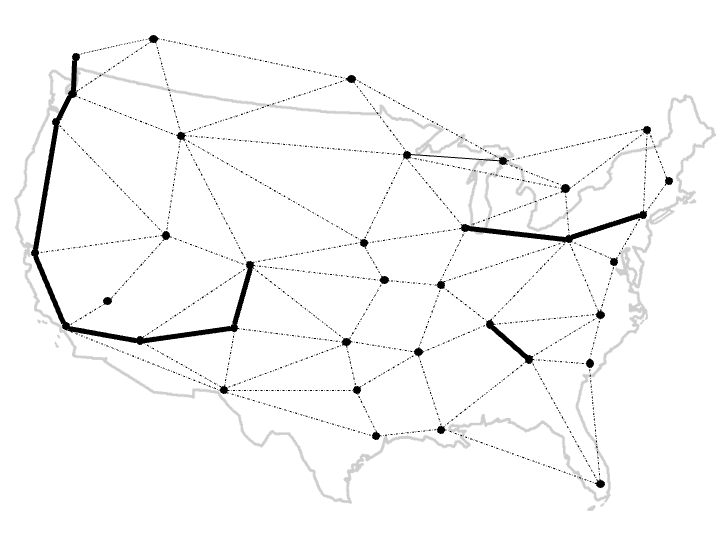}
\caption{The ten most frequently chosen routes.} 
\label{FrequentlyChosenRoutesMap}
\end{figure}

Figure \ref{FrequentlyChosenRoutesMap} shows the ten most frequently selected routes. Interestingly, all ten of these routes are chosen in the optimal solution with 45 train cars. Seven of these ten most frequently selected routes are entirely west of the Rocky Mountains, which may reflect the relative sparsity of routes in the west, as well as their importance in satisfying high-value tickets. 
\begin{table}[]
    \centering 
    \begin{tabular}{|l|c|c|}
    \hline
    Ticket & Points & Frequency Chosen \\
    \hline
    Atlanta to New York & 6 & 25 \\
    Chicago to Santa Fe & 9 & 28 \\
    Los Angeles to Seattle & 9 & 21 \\
    Phoenix to Portland & 11 & 23 \\
    Santa Fe to Vancouver & 13 & 21 \\
    Chicago to Los Angeles & 16 & 28 \\
    Nashville to Portland & 17 & 21 \\
    Atlanta to San Francisco & 17 & 21 \\
    New York to Los Angeles & 21 & 30 \\
    New York to Seattle & 22 & 21 \\
    \hline
 \end{tabular}
    \caption{\color{black} The ten most frequently chosen tickets (out of 50 total instances).} \label{CardFrequency}
\end{table}

\section{Conclusion}
We \textcolor{black}{model the maximum} score in the popular board game Ticket to Ride using an undirected graph,  a directed graph, and then as a mixed-integer program. Our numerical study  determines the optimal score parameterized by the number of train cars.  We describe the most frequently chosen tickets and routes within these 50 instances. Ticket to Ride appears amenable to further mathematical analyses, including its probabilistic and competitive aspects, which are beyond the scope of this paper.

\bibliographystyle{ieeetr}
\bibliography{bibliography}
\end{document}